\documentclass[a4paper,12pt]{amsart}
\usepackage{amssymb,amsmath,amsthm}
\usepackage{hyperref}
\usepackage[utf8]{inputenc}
\usepackage[T1]{fontenc}
\usepackage{microtype}
\usepackage[english,polish]{babel}
\newtheorem{theorem}{Theorem}[section]
\newtheorem{lemma}[theorem]{Lemma}
\newtheorem{corollary}[theorem]{Corollary}
\theoremstyle{definition}
\newtheorem{remark}[theorem]{Remark}
\newtheorem{example}[theorem]{Example}
\begin{document}

\selectlanguage{polish}
\author[J. Morawiec]{Janusz Morawiec}
\address{Instytut Matematyki{}\\
Uniwersytet Śląski{}\\
Bankowa 14, PL-40-007 Katowice{}\\
Poland}
\email{morawiec@math.us.edu.pl}

\author[T. Zürcher]{Thomas Zürcher}
\address{Instytut Matematyki{}\\
Uniwersytet Śląski{}\\
Bankowa 14, PL-40-007 Katowice{}\\
Poland}
\email{thomas.zurcher@us.edu.pl}
\selectlanguage{English}
	
\title[Functional equations for generating $\varepsilon$-invariant measures]{An application of functional equations for generating $\varepsilon$-invariant measures}
\begin{abstract}
Let $(X,{\mathcal A},\mu)$ be a probability space and let \mbox{$S\colon X\to X$} be a measurable transformation. Motivated by the paper of K.~Nikodem [Czechoslovak Math. J. 41(116) (4) (1991) 565--569], we concentrate on a functional equation generating measures that are absolutely continuous with respect to $\mu$ and $\varepsilon$\nobreakdash-invariant under $S$. As a consequence of the investigation, we obtain a result on the existence and uniqueness of solutions $\varphi\in L^1([0,1])$ of the functional equation
$$
\varphi(x)=\sum_{n=1}^{N}|f_n'(x)|\varphi(f_n(x))+g(x),
$$
where $g\in L^1([0,1])$ and $f_1,\ldots,f_N\colon[0,1]\to[0,1]$ are functions satisfying some extra conditions.
\end{abstract}
\maketitle

%%%%%%%%%%%%%%%% Section 1 %%%%%%%%%%%%%%%%%%

\renewcommand{\theequation}{1.\arabic{equation}}\setcounter{equation}{0}
\section{Introduction}
The aim of this paper is to study the problem of the existence of solutions $\varphi\in L^1(X)$ of the following equation
\begin{equation}\label{E}
\varphi=P\varphi+g,
\end{equation}
where $g\in L^1(X)$ and $P\colon L^1(X)\to L^1(X)$ are given.
In the rest of the introduction, we let $X=[0,1]$ be equipped with the Borel $\sigma$\nobreakdash-algebra and the Lebesgue measure.
The motivation for studying such a problem is twofold. An original impulse for our investigation came from the
paper \cite{N1991}, where integrable solutions of the equation
\begin{equation}\label{e0}
\varphi(x)=\frac{1}{2}\varphi\left(\frac{x}{2}\right)+\frac{1}{2}\varphi\left(\frac{x+1}{2}\right)+g(x)
\end{equation}
were investigated in connection with $\varepsilon$\nobreakdash-invariant measures under the $2$\nobreakdash-adic transformation. The next section contains more details on $\varepsilon$\nobreakdash-\hspace{0pt}invariant measures and on functional equations associated with them. Let us note here only that equation \eqref{e0} is a very particular case of an interesting functional equation of the form
\begin{equation}\label{e}
\varphi(x)=\sum_{n=1}^{N}|f_n'(x)|\varphi(f_n(x))+g(x).
\end{equation}
We will always assume that $N\geq 2$.

The second inspiration to study integrable solutions of equation \eqref{e}, and hence also of \eqref{E}, is strictly connected with a problem posed by Janusz Matkowski in \cite{M1985} and a question posed during the 47th International Symposium on Functional Equations by Jacek Wesołowski in connection with probability measures investigated in \cite{MW2012}. Namely, assume that $f_1,\ldots,f_N\colon[0,1]\to[0,1]$ are strictly increasing contractions satisfying
\begin{equation}\label{c2}
f_n((0,1))\cap f_m((0,1))=\emptyset\quad\hbox{for all }n\neq m
\end{equation}
and consider the class $\mathcal C$ consisting of all increasing and continuous functions $\phi\colon [0,1]\to [0,1]$ such that $\phi(0)=0$, $\phi(1)=1$ and
\begin{equation*}
\phi(x)=\sum_{n=1}^{N}\phi(f_n(x))-\sum_{n=1}^{N}\phi(f_n(0)).
\end{equation*}
Wide classes of singular functions belonging to the class $\mathcal C$ were constructed in \cite{MZ2018,MZ}. So far we have some idea how singular functions from the class $\mathcal C$ look like, but we do not know much about absolutely continuous functions from $\mathcal C$. However, we know that under a weak assumption on the functions $f_n$, each function belonging to the class $\mathcal C$ can be expressed as a convex combination of absolutely continuous and singular functions from $\mathcal C$. Finally, observe that (again under weak assumptions on the functions $f_n$) absolutely continuous functions belonging to the class $\mathcal C$ are in one-to-one correspondence with densities satisfying \eqref{e} with~$g=0$.

%%%%%%%%%%%%%%%% Section 2 %%%%%%%%%%%%%%%%%%

\renewcommand{\theequation}{2.\arabic{equation}}\setcounter{equation}{0}
\section{Preliminaries}
In this section we explain more precisely our impulse for studying integrable solutions of equation \eqref{e}, as well as of its generalization~\eqref{E}. We begin with recalling some definitions and results useful in the main part of this paper.

Throughout this paper we assume that $(X,{\mathcal A},\mu)$ is a probability space and $S\colon X\to X$ is a measurable transformation.

In the case where $X$ is a Borel subset of $\mathbb R$ we assume that $\mathcal{A}$ is the $\sigma$-algebra $\mathcal B(X)$ of all Borel subsets of $X$ and $\mu$ is the Lebesgue measure restricted to $\mathcal B(X)$.

We say that $S$ is \emph{nonsingular} if $\mu(S^{-1}(A))=0$ for every $A\in\mathcal A$ such that $\mu(A)=0$. We say that $S$ is \emph{measure preserving} if $\mu(S^{-1}(A))=\mu(A)$ for every $A\in\mathcal A$; we will alternately say that the measure $\mu$ is \emph{invariant} under $S$ if $S$ is measure preserving. Observe that every measure preserving transformation is nonsingular.

Fix a real number $\varepsilon\geq 0$. A probability measure $\nu$ defined on ${\mathcal A}$ is said to be \emph{$\varepsilon$-invariant under $S$} if
\begin{equation*}
|\nu(S^{-1}(A))-\nu(A)|\leq\varepsilon\mu(A)
\end{equation*}
for every $A\in{\mathcal A}$. It is clear that every measure that is $0$-invariant under~$S$ is invariant under $S$, and so the concept of measures $\varepsilon$-invariant under~$S$ generalizes the notion of measures invariant under $S$.

From now on, for a nonsingular $S$ we denote by $P_S$ the corresponding \emph{Frobenius--Perron operator}, i.e.\ $P_S\colon L^1(X)\to L^1(X)$ is the operator uniquely defined by the equation
\begin{equation}\label{int}
\int_A P_Sf(x)d\mu(x)=\int_{S^{-1}(A)}f(x)d\mu(x)\quad\hbox{for every }A\in{\mathcal A}.
\end{equation}
The operator $P_S$ is linear and continuous. If $S$ is nonsingular, then every $m$-th iterate $S^m$ of $S$ is also nonsingular and the Frobenius-Perron operator corresponding to $S^m$ is the $m$-th iterate $P_S^m$ of $P_S$. Here and throughout we adopt the convention that $P^0={\rm id}_{X}$ for every operator $P\colon L^1(X)\to L^1(X)$. In the case where $X=[0,1]$ and $\mu$ is the one-dimensional Lebesgue measure the Frobenius--Perron operator corresponding to $S$ can be written explicitly as follows
\begin{equation}\label{explicit}
P_Sf(x)=\frac{d}{dx}\int_{S^{-1}([0,x])}f(y)\,dy
\end{equation}
(see \cite[Formula 1.2.7]{LM1994}).

A useful tool for studying measures $\varepsilon$-invariant under nonsingular $S$ that are absolutely continuous with respect to $\mu$ reads as follows.
\begin{theorem}[see {\cite[Theorem~1]{N1991}}]\label{thm21}
A finite measure $\nu$ that is absolutely continuous with respect to $\mu$ is $\varepsilon$-invariant under a nonsingular $S$ if and only if the Radon--Nikodym derivative $f$ of $\nu$ with respect to $\mu$ satisfies
\begin{equation*}
|P_Sf(x)-f(x)|\leq\varepsilon\quad\hbox{for $\mu$-almost all }x\in X.
\end{equation*}
\end{theorem}
In the case where $\varepsilon=0$ in Theorem~\ref{thm21}, we recognize the well know fact saying that an absolutely continuous measure $\nu$ with respect to $\mu$ is invariant under a nonsingular $S$ if and only if the Radon-Nikodym derivative $f$ of $\nu$ with respect to $\mu$ is a fixed point of the Frobenius--Perron operator $P_S$ (see \cite[Theorem~4.1.1]{LM1994}).

From Theorem~\ref{thm21} we see that to find measures $\varepsilon$-invariant under a nonsingular transformation $S$ it is enough to solve, in the space $L^1(X)$, equation~\eqref{E} with $P=P_S$ and $|g|\leq\varepsilon$. Therefore, we are interested in finding solutions $\varphi\in L^1(X)$ of the following special case of equation~\eqref{E}
\begin{equation}\label{PE}
\varphi=P_S\varphi+g.
\end{equation}
We will do it for a wide class of important transformations. But before we introduce the class, observe that according to~\eqref{int} a necessary condition for $g\in L^1(X)$ in order that equation~\eqref{PE} has a solution $\varphi\in L^1(X)$ is
\begin{equation}\label{g}
\int_Xg(x)d\mu(x)=0.
\end{equation}

A measure preserving $S$ such that $S(A)\in\mathcal A$ for every $A\in\mathcal A$ is said to be \emph{exact} if $\lim_{m\to\infty}\mu\big(S^m(A)\big)=1$ for every $A\in\mathcal A$ such that $\mu(A)>0$. The following characterization of exactness is well known.

\begin{theorem}[see {\cite[Corollary~4.4.1]{LM1994}}]\label{thm22}
A measure preserving $S$ is exact if and only if for every $f\in L^1(X)$ the sequence $(P_S^mf)_{m\in\mathbb N}$ converges in~$L^1(X)$ to $\int_Xf(x)d\mu(x)$.
\end{theorem}

%%%%%%%%%%%%%%%% Section 3 %%%%%%%%%%%%%%%%%%

\renewcommand{\theequation}{3.\arabic{equation}}\setcounter{equation}{0}
\section{A generalization of a Nikodem result}
We begin this section with a result on the existence of integrable solutions of equation~\eqref{PE}, whose proof is a direct trace of the proof of Theorem~2 from \cite{N1991}, but for the convenience of the readers we repeat it.

\begin{theorem}\label{thm31}
Assume that $g\in L^1(X)$, $S$ is exact and $P_S1=1$. Then equation $\eqref{PE}$ has a solution in
$L^1(X)$ if and only if the series $\sum_{m=0}^{\infty}P_S^m g$ converges in $L^1(X)$. Moreover, every solution $\varphi\in L^1(X)$ of equation $\eqref{PE}$ is of the form
\begin{equation*}
\varphi=\sum_{m=0}^{\infty}P_S^m g+c,
\end{equation*}
where $c$ is a real constant.
\end{theorem}

\begin{pf}	
Assume first that the series $\sum_{m=0}^{\infty}P_S^m g$ converges in $L^1(X)$. Fix a real constant $c$ and set $\varphi=\sum_{m=0}^{\infty}P_S^m g+c$. The linearity and continuity of $P_S$ jointly with the equality $P_S1=1$ imply
\begin{equation*}
P_S\varphi+g=\sum_{m=0}^{\infty}P_S^{m+1}g+P_Sc+g=\sum_{m=0}^{\infty}P_S^m g+c=\varphi.
\end{equation*}
	
Assume now that $\varphi\in L^1(X)$ satisfies~\eqref{PE}. Then, by the linearity of~$P_S$, we have
\begin{equation*}
P_S^k\varphi=P_S^{k+1}\varphi+P_S^kg
\end{equation*}
for every $k\in\mathbb N$. Adding the above equation over $k=0,\ldots,m$ leads to
\begin{equation*}
\sum_{k=0}^m P_S^k g=\varphi-P_S^{m+1}\varphi.
\end{equation*}
for every $m\in\mathbb N$. Finally, passing with $m$ to $\infty$ and making use of Theorem~\ref{thm22}, we conclude that the series $\sum_{k=0}^\infty P_S^k g$ converges in $L^1(X)$ and that
\begin{equation*}
\sum_{k=0}^\infty P_S^k g=\varphi-\int_X\varphi(x)d\mu(x),
\end{equation*}
which completes the proof.\hfill$\square$	
\end{pf}

Now we define a transformation $S$ whose Frobenius--Perron operator appears in equation~\eqref{e}. For this purpose  put $X=[0,1]$ and fix strictly monotone functions $f_1,\ldots,f_N\colon[0,1]\to[0,1]$ satisfying~\eqref{c2}. Note that the functions $f_n$ are differentiable almost everywhere on $[0,1]$ as they are monotone. Now define the announced transformation $S\colon[0,1]\to[0,1]$ by putting
\begin{equation}\label{S}
S(x)=\begin{cases}
f_n^{-1}(x)&\hbox{ for }x\in f_n((0,1))\hbox{ and }n\in\{1,\ldots,N\},\\
0&\hbox{ for }x\not\in\bigcup_{n=1}^{N}f_n((0,1)).
\end{cases}
\end{equation}
Clearly, $S$ is well defined by~\eqref{c2}. Moreover, it is nonsingular provided that all the functions $f_1,\ldots,f_N$ satisfy Luzin's condition~N (or equivalently all the inverses $f_1^{-1},\ldots,f_N^{-1}$ are nonsingular) and
\begin{equation}\label{c1}
\bigcup_{n=1}^{N}f_n([0,1])=[0,1].
\end{equation}
Applying formula~\eqref{explicit} we see that the Frobenius--Perron operator $P_S$ corresponding to a nonsingular $S$ defined by formula~\eqref{S} is of the form
\begin{equation}\label{PS}
P_S f=\sum_{n=1}^{N}|f_n'|\big(f\circ f_n\big).
\end{equation}
Therefore, in the case where $S$ is defined by formula~\eqref{S} equation~\eqref{PE} reduces to equation~\eqref{e} and Theorem \ref{thm31} implies the following result.

\begin{corollary}\label{cor32}
Assume that $g\in L^1([0,1])$ and $f_1,\ldots,f_N\colon[0,1]\to[0,1]$ are strictly monotone nonsingular functions satisfying~\eqref{c2} and
\begin{equation}\label{f'}
\sum_{n=1}^{N}|f_n'(x)|=1\quad\hbox{ for almost all }x\in[0,1].
\end{equation}
If $S\colon[0,1]\to[0,1]$ defined by formula~\eqref{S} is exact, then equation~\eqref{e} has a solution in $L^1([0,1])$ if and only if the series
\begin{equation*}\label{formula}
\sum_{m=1}^{\infty}\sum_{n_1,\ldots,n_m=1}^N\left(\prod_{k=1}^{m}|f_{n_k}'\circ f_{n_{k-1}}\circ\cdots\circ f_{n_1}|\right)g\circ f_{n_m}\circ f_{n_{m-1}}\circ\cdots\circ f_{n_1}
\end{equation*}
converges in $L^1([0,1])$. Moreover, every solution $\varphi\in L^1([0,1])$ of equation~\eqref{e} is of the form
\begin{equation*}
\varphi=g+\!\sum_{m=1}^{\infty}\sum_{n_1,\ldots,n_m=1}^N\!\!\left(\prod_{k=1}^{m}|f_{n_k}'\!\circ\! f_{n_{k-1}}\!\circ\cdots\circ\! f_{n_1}|\!\!\right)\!g\circ f_{n_m}\circ f_{n_{m-1}}\circ\cdots\circ f_{n_1}+c,
\end{equation*}
where $c$ is a real constant.
\end{corollary}

There is a known criterion (being very close to a characterization) of exactness of nonsingular transformations (see \cite[Proposition~5.6.2 and Remark~5.6.1 on p.~111\footnote{There are actually two remarks numbered 5.6.1 in the book.}]{LM1994}). A wide class of interesting examples of exact transformations defined by formula~\eqref{S}, including transformations studied by R\'{e}nyi in \cite{R1957} and by Rohlin in \cite{Ro1961} (see also \cite{rohlin1964exact}), can be found in \cite[Chapter 6.2]{LM1994}. Now we give one of the simplest possible realization of the assumptions of Corollary~\ref{cor32}. For this purpose, fix non-zero real numbers $\alpha_1,\ldots,\alpha_N$ such that
\begin{equation}\label{1}
\sum_{n=1}^N|\alpha_n|=1
\end{equation}
and put
\begin{equation*}
\beta_n=\sum_{k=1}^{n}|\alpha_k|-\frac{|\alpha_{n}|+\alpha_{n}}{2},
\end{equation*}
\begin{equation}\label{affine}
f_n(x)=\alpha_n x+\beta_n
\end{equation}
for all $x\in[0,1]$ and $n\in\{1,\ldots,N\}$. Clearly, \eqref{1} implies~\eqref{f'}. Moreover, it is well known that in the considered case the transformation defined by formula~\eqref{S} is exact (see \cite[Theorem~6.2.1, Definition~5.6.2 and Proposition~5.6.2]{LM1994}). Thus by Corollary~\ref{cor32} we conclude that any solution $\varphi\in L^1([0,1])$ of the equation
\begin{equation}\label{e1}
\varphi(x)=\sum_{n=1}^{N}|\alpha_n|\varphi\left(\alpha_n x+\beta_n\right)+g(x)
\end{equation}
is of the form
\begin{align}\label{sol}
\varphi(x)=&\sum_{m=1}^{\infty}\sum_{n_1,\ldots,n_m=1}^N\left(\prod_{k=1}^{m}|\alpha_{n_k}|\right)g\left(\prod_{i=1}^{m}\alpha_{n_i}x+\sum_{i=1}^{m}\beta_{n_i}\prod_{j=i+1}^{m}\alpha_{n_j}\right)\nonumber\\
&+g(x)+c,
\end{align}
where $c\in\mathbb R$; here we adopt the convention that $\prod_{j=m+1}^{m}a_{j}=1$ for all $m\in\mathbb N$ and $a_j\in\mathbb R$.

Finally, note that Corollary \ref{cor32} reduces to Theorem~2 from \cite{N1991} in the case where $N=2$ and $f_1,f_2$ are given by~\eqref{affine} with $\alpha_1=\alpha_2=\frac{1}{2}$.

%%%%%%%%%%%%%%%% Section 4 %%%%%%%%%%%%%%%%%%

\renewcommand{\theequation}{4.\arabic{equation}}\setcounter{equation}{0}
\section{Examples}
In this section we give four examples of measures $\varepsilon$-invariant under a given nonsingular transformation. The first two examples will be constructed by hand, whereas in the next two we will apply Theorem~\ref{thm31} and Corollary \ref{cor32}.

We begin with the general observation that every measure invariant under a nonsingular transformation $S$ generates a large class of measures that are $\varepsilon$-invariant under $S$.

\begin{example}\label{ex1}
Assume that $S$ is nonsingular and $\mu$ is invariant under~$S$. Fix sets $A_1,\ldots,A_m\in\mathcal A$ and real numbers  $\varepsilon, \varepsilon_1,\ldots,\varepsilon_m>0$ such that $\varepsilon=\sum_{i=1}^m\varepsilon_i$. Next define a finite measure $\nu$ on $\mathcal A$ by putting
\begin{equation*}
\nu(A)=\sum_{i=1}^m\varepsilon_i\mu(A\cap A_i).
\end{equation*}
Then observe that for every $A\in{\mathcal A}$ we have
\begin{align*}
\left|\nu\big(S^{-1}(A)\big)-\nu(A)\right|&\leq\sum_{i=1}^m\varepsilon_i\left|\mu\big(S^{-1}(A)\cap A_i\big)-\mu(A\cap A_i)\right|\\
&\leq\sum_{i=1}^m\varepsilon_i\max\left\{\mu\big(S^{-1}(A)\cap A_i\big),\mu(A\cap A_i)\right\}\\
&\leq\sum_{i=1}^m\varepsilon_i\max\left\{\mu\big(S^{-1}(A)\big),\mu(A)\right\}=\varepsilon\mu(A).
\end{align*}
\end{example}

Note that if $x_1,\ldots,x_m$ are fixed points of a nonsingular $S$, then every convex combination of the Dirac measures $\delta_{x_1},\ldots,\delta_{x_m}$ is a probability measure that is invariant under $S$. Thus Example \ref{ex1} shows that any transformation defined by formula~\eqref{S} with strictly monotone functions $f_1,\ldots,f_N$ satisfying~\eqref{c2}, \eqref{c1} and sending sets of measure zero to sets of measure zero has $\varepsilon$-invariant measures.

The next example is of similar type as the first one, however the difference is that we will not assume the existence of a measure that is invariant under~$S$. It concerns transformations defined by formula~\eqref{S} with $X=[0,1)$ and $f_n$ of the form~\eqref{affine}.

\begin{example}\label{ex2}
Fix positive real numbers $\alpha_1,\ldots,\alpha_N$ satisfying~\eqref{1} and let for every $n\in\{1,\ldots,N\}$ the function $f_n$ be given by~\eqref{affine}. Then put $X=[0,1)$ and consider the transformation $S\colon [0,1)\to[0,1)$ defined by~\eqref{S}. Obviously, $S$ is nonsingular.
	
For all $k\in\mathbb N$ and $n_1,\ldots,n_k\in\{1,\ldots,N\}$ we put
\begin{equation*}
I_{n_k,\ldots,n_1}=[f_{n_k}\circ\cdots\circ f_{n_1}(0),f_{n_k}\circ\cdots\circ f_{n_1}(1)).
\end{equation*}
A simple induction with the use of the strict increasingness of $f_1,\ldots,f_N$, \eqref{c2}~and~\eqref{c1} shows that for every $k\in\mathbb N$ and $n_1,\ldots,n_k,m_1,\ldots,m_k\in\{1,\ldots,N\}$ we have
\begin{equation}\label{w3}
I_{n_k,\ldots,n_1}\cap I_{m_k,\ldots,m_1}=\emptyset\quad\hbox{ for }(n_k,\ldots,n_1)\neq (m_k,\ldots,m_1),
\end{equation}
\begin{equation}\label{w1}
I_{n_k,\ldots,n_1}=\bigcup_{n=1}^{N}I_{n_k,\ldots,n_1,n},
\end{equation}
\begin{equation}\label{w2}
\bigcup_{n_1,\ldots,n_k=1}^{N}I_{n_k,\ldots,n_1}=[0,1).
\end{equation}
	
Next for every  $k\in\mathbb N$ we put
\begin{equation*}
\mathcal S_k=\big\{I_{n_k,\ldots,n_1}:n_1,\ldots,n_k\in\{1,\ldots,N\}\big\}\quad\hbox{and}\quad \mathcal  S_0=\{[0,1)\}.
\end{equation*}
It is easy to see that the family $\mathcal S=\bigcup_{k=0}^\infty\mathcal S_k\cup\{\emptyset\}$ is a semi-algebra of subsets of the interval $[0,1)$; i.e. $[0,1)\in\mathcal S$ and $I,J\in\mathcal S$ implies $I\cap J\in\mathcal S$ and $[0,1)\setminus I$ can be expressed as a finite disjoint union of sets in $\mathcal S$ (see \cite[Definition~1.4.1]{P1977}). Note that the $\sigma$-algebra generated by the semi-algebra $\mathcal S$ coincides with the family ${\mathcal B}([0,1))$ of all Borel subsets of the interval $[0,1)$.
	
Fix $\varepsilon\in[0,1]$, two different numbers $p,q\in\{1,\ldots,N\}$ and define a function $\xi\colon\{0,\ldots,N\}\to\mathbb R$ by putting
\begin{equation*}
\xi(n)=\begin{cases}
\varepsilon\min\{\alpha_p,\alpha_q\}&\hbox{ for }n=p,\\
-\varepsilon\min\{\alpha_p,\alpha_q\}&\hbox{ for }n=q,\\
0&\hbox{ for }n\not\in\{p,q\}.
\end{cases}
\end{equation*}
Clearly,
\begin{equation}\label{xi}
\sum_{n=1}^{N}\xi(n)=0.
\end{equation}
Now we want to define a probability measure $\nu_0$ on $\mathcal S$; i.e.\ a function $\nu_0\colon\mathcal S\to[0,1]$ such that $\nu_0(\emptyset)=0$ and $\nu_0(\bigcup_{n\in\mathbb N}J_n)=\sum_{n\in\mathbb N}\nu_0(J_n)$ for all pairwise disjoint elements $(J_n)_{n\in\mathbb N}$ of $\mathcal S$ with $\bigcup_{n\in\mathbb N}J_n\in\mathcal S$ (see \cite[Section 2.3]{P1977}). We will do it inductively.
	
In the first step we define $\nu_0$ on $\mathcal  S_0\cup\mathcal S_1$ by putting
\begin{equation*}
\nu_0([0,1))=1\quad \quad\hbox{ and }\quad\nu_0(I_n)=\alpha_n+\xi(n).
\end{equation*}
It is clear that $\nu_0(I_n)\in[0,1]$ for every $I_n\in\mathcal  S_1$. Moreover, \eqref{1}~and~\eqref{xi} imply
\begin{equation*}\label{step1}
\nu_0\left([0,1)\right)=1=\sum_{n=1}^{N}\alpha_n=\sum_{n=1}^{N}\nu_0\left(I_{n}\right),
\end{equation*}
and we see (according to \eqref{w2}~and~\eqref{w3} with $k=1$) that $\nu_0$ is well defined on $\mathcal S_0\cup\mathcal  S_1$.
	
To do the inductive step we assume that for a fixed $k\in\mathbb N$ we have defined $\nu_0$ on $\bigcup_{n=0}^{k}\mathcal  S_k$. Now we define $\nu_0$ on $\mathcal  S_{k+1}$ by putting
\begin{equation*}
\nu_0(I_{n_{k+1},\ldots,n_1})=\big(\alpha_{n_{k+1}}+\xi(n_{k+1})\big)\prod_{i=1}^{k}\alpha_{n_i}.
\end{equation*}
(Note that adopting the convention that $\prod_{i=1}^{0}\alpha_{n_i}=1$ the above formula for $\nu_0$ coincides with that from the first step.) Applying~\eqref{1} we obtain
\begin{equation}\label{suma}
\sum_{n=1}^{N}\nu_0(I_{n_{k},\ldots,n_1,n})=\big(\alpha_{n_k}+\xi(n_k)\big)\prod_{i=1}^{k-1}\alpha_{n_i}\sum_{n=1}^{N}\alpha_n=\nu_0(I_{n_{k},\ldots,n_1})
\end{equation}
for every $I_{n_{k},\ldots,n_1}\in\mathcal S_k$, which means (according to \eqref{w1}~and~\eqref{w3}) that $\nu_0$ is well defined on $\bigcup_{n=0}^{k+1}\mathcal  S_k$.
	
Finally, putting $\nu_0(\emptyset)=0$ we have defined $\nu_0$ on $\mathcal S$.
	
It remains to prove that $\nu_0$ is $\sigma$-additive. For this purpose fix a pairwise disjoint sequence $(J_m)_{m\in\mathbb N}$ of elements of the semi-algebra $\mathcal S$ such that $\bigcup_{m\in\mathbb N}J_m=J\in\mathcal S$. It simplifies the argument, and causes no loss of generality, to assume $J=[0,1)$. Thus we need to show that $\sum_{m\in\mathbb N}\nu_0(J_m)=1$.
	
Define a nondecreasing sequence $(k_m)_{m\in\mathbb N}$ of integers in such a way that $\{J_1,\ldots,J_m\}\subset\bigcup_{k=1}^{k_m}\mathcal S_k$ for every $m\in\mathbb N$. Note that for all $m\in\mathbb N$,  $l\in\{1,\ldots,m\}$ and $I\in \mathcal S_{k_m}$, we have either $J_l\cap I=\emptyset$ or $J_l\cap I=I$. Next for every $m\in\mathbb N$ put
\begin{equation*}
\mathcal D_m=\left\{I\in \mathcal S_{k_m}: I\subset\bigcup_{l=1}^{m}J_l\right\}\quad\hbox{ and }\quad
d_m=\sum_{I\in\mathcal S_{k_m}\setminus\mathcal D_m}\nu_0(I).
\end{equation*}
Making use of \eqref{w1}, \eqref{suma} and \eqref{w2} we conclude that
\begin{equation*}
\sum_{l=1}^m\nu_0(J_l)=\sum_{I\in \mathcal D_m}\nu_0(I)=\sum_{I\in\mathcal S_{k_m}}\nu_0(I)-d_m=1-d_m.
\end{equation*}
Now it is enough to show that
\begin{equation}\label{dm}
\lim_{m\to\infty}d_m=0.
\end{equation}
	
By the definition of $\nu_0$ it is easy to see that for every $I\in\mathcal S$ we have
\begin{equation*}
\nu_0(I)\leq 2l(I);
\end{equation*}
here and later on the symbol $l$ denotes the Lebesgue measure on the real line. This jointly with~\eqref{w3} yields
\begin{equation*}
d_m\leq 2\sum_{I\in\mathcal S_{k_m}\setminus\mathcal D_m}l(I)=2l\left(\bigcup_{I\in\mathcal S_{k_m}\setminus\mathcal D_m}I\right)=2l\left(\bigcup_{l=m+1}^{\infty} J_l\right).
\end{equation*}
Passing with $m$ to $\infty$ we get~\eqref{dm}.
	
Thus we have proved that $\nu_0$ is a probability measure defined on the semi-algebra $\mathcal S$.
	
Extend $\nu_0$ to a probability measure $\nu\colon\mathcal B([0,1))\to[0,1]$; such an extension exists and it is unique (see \cite[Corollary~2.4.9 and Proposition~2.5.1]{P1977}).
	
Now we will show that the measure $\nu$ is $\varepsilon$-invariant under $S$.
	
For this purpose fix $I_{n_k,\ldots,n_1}\in \mathcal{S}_k$. Then
\begin{equation*}
S^{-1}(I_{n_k,\ldots,n_1})=\bigcup_{n=1}^{N}I_{n,n_k,\ldots,n_1},
\end{equation*}
which jointly with \eqref{w3}, \eqref{1} and \eqref{xi} implies
\begin{align*}
\nu\big(S^{-1}(I_{n_k,\ldots,n_1})\big)&=\sum_{n=1}^{N}\nu_0(I_{n,n_k,\ldots,n_1})=\prod_{i=1}^{k}\alpha_{n_i}\sum_{n=1}^{N}\big(\alpha_{n}+\xi(n)\big)\\
&=\prod_{i=1}^{k}\alpha_{n_i}=\nu(I_{n_k,\ldots,n_1})-\xi(n_k)\prod_{i=1}^{k-1}\alpha_{n_i}.
\end{align*}
In consequence,
\begin{equation}\label{nuS}
\left|\nu\big(S^{-1}(I_{n_k,\ldots,n_1})\big)-\nu(I_{n_k,\ldots,n_1})\right|=|\xi(n_k)|\prod_{i=1}^{k-1}\alpha_{n_i}\leq\varepsilon l(I_{n_{k},\ldots,n_1}).
\end{equation}
	
Fix a set $B\in{\mathcal B}([0,1))$, a number $\delta>0$ and choose a countable family $\{F_j:j\in J\}$ of pairwise disjoint elements of the semi-algebra $\mathcal S$ such that
$\bigcup_{j\in J}F_j\subset B$,
\begin{equation*}
\left|\nu(B)-\nu\left(\bigcup_{j\in J}F_j\right)\right|\leq\delta
\;\:\hbox{and}\;\:
\left|\nu(S^{-1}(B))-\nu\left(S^{-1}\left(\bigcup_{j\in J}F_j\right)\right)\right|\leq\delta;
\end{equation*}
such a family exists, because on any complete separable metric space any finite Borel measure is regular (see \cite[Theorem~7.1.4]{Dud2002}). Then making use of~\eqref{nuS} we get
\begin{align*}
\left|\nu\big(S^{-1}(B)\big)-\nu(B)\right|&\leq
\sum_{j\in J}\left|\nu\left(S^{-1}(F_j)\right)-\nu(F_j)\right|+2\delta\\
&\leq\varepsilon\sum_{j\in J} l(F_j)+2\delta
\leq\varepsilon l(B)+2\delta.
\end{align*}
Finally, tending with $\delta$ to $0$ we conclude that the measure $\nu$ is $\varepsilon$-invariant under $S$.
\end{example}

To give the next two examples we need the following observation.

\begin{remark}\label{rem43}$ $
\begin{enumerate}
\item If $\varphi$ solves~\eqref{PE}, then $\int_{X}g(x)d\mu(x)=0$.
\item If $P_Sg=0$, then $\int_{X}g(x)d\mu(x)=0$.
\item If $P_Sg=0$, then $P_S^mg=0$ for every $m\in\mathbb N$.
\end{enumerate}
\end{remark}

\begin{pf}
The first two statements are an immediate consequence of~\eqref{int}, whereas the third one follows from the linearity of $P_S$.\hfill$\square$
\end{pf}

Remark~\ref{rem43} helps us to apply Theorem~\ref{thm31}. Indeed, assertion (i) says that to find an integrable solution of~\eqref{PE} we must assume that $g$ has integral equals zero over $X$. In view of assertion (ii) it can be realized by choosing $g$ in such a way that $P_Sg=0$. Having chosen such a $g$, assertion~(iii) implies the convergence of the series $\sum_{m=0}^{\infty}P_S^m g$ to~$g$. Concluding, if we fix $\varepsilon\in[0,1]$ and $g\in L^1(X)$ such that $|g|\leq\varepsilon$ and $P_Sg=0$, then $g+1$ is the density of a probability measure that is $\varepsilon$-invariant under $S$, i.e.\ the formula
\begin{equation}\label{nu}
\nu(A)=\mu(A)+\int_{A}g(x)d\mu(x)\quad\hbox{for every }A\in\mathcal A
\end{equation}
defines a probability measure that is $\varepsilon$-invariant under $S$.

To see that in many cases $g$ can be fixed in such a way that $|g|\leq\varepsilon$ and $P_Sg=0$, we consider in the next two examples the case where $X=[0,1]$ and $S$ is defined by formula~\eqref{S}.

\begin{example}\label{ex3}
Fix strictly monotone functions $f_1,\ldots,f_N\colon[0,1]\to[0,1]$ satisfying \eqref{c2}, \eqref{c1}, \eqref{f'}.
We assume that $f_N(0)\geq f_n(1)$ for $n\leq N-1$ and further that $|f_N'|\geq\frac{1}{2}$. Then choose an integrable function $g_0\colon [0,f_{N}(0)]\to[0,\varepsilon]$ and extend it to an integrable function $g\colon[0,1]\to[0,\varepsilon]$ by putting $g=-\frac{1}{|f_N'|}\sum_{n=1}^{N-1}|f_n'|\big(g_0\circ f_n\circ f_N^{-1}\big)$ on $(f_{N}(0),1]$. Since in the considered case $P_S$ is of the form~\eqref{PS}, we have $P_Sg=\sum_{n=1}^{N}|f_n'|\big(g\circ f_n\big)=0$.
\end{example}

The last example shows not only how to choose a $g$ with $|g|\leq\varepsilon$ and $P_Sg=0$, but how to choose such a $g$ to calculate the integral in~\eqref{nu}.

\begin{example}\label{ex4}	
Put $X=[0,1]$ and let $S$ be defined by formula~\eqref{S} with the $f_n$ given by~\eqref{affine}, where $\alpha_1,\ldots,\alpha_N$ are non-zero real numbers satisfying~\eqref{1}. Fix $\varepsilon\in[0,1]$ and real numbers $\gamma_1,\ldots,\gamma_{N}\in[-\varepsilon,\varepsilon]$ such that $\sum_{n=1}^N |\alpha_n|\gamma_n=0$. Choose $g$ to be a constant equal to $\gamma_n$ on every interval $I_n=(\min\{f_n(0),f_n(1)\},\max\{f_n(0),f_n(1)\})$. Then $P_Sg=\sum_{n=1}^N|\alpha_n|\gamma_n=0$. Finally, according to~\eqref{nu} we conclude that the formula
\begin{equation*}
\nu(A)=\sum_{n=1}^N(1+\gamma_n)l(A\cap I_n)
\end{equation*}
defines a Borel probability measure that is $\varepsilon$-invariant under $S$.
\end{example}

%%%%%%%%%%%%%%%% Section 5 %%%%%%%%%%%%%%%%%%

\renewcommand{\theequation}{5.\arabic{equation}}\setcounter{equation}{0}
\section{Further results}
We begin this section with a generalization of Theorem~\ref{thm31}. To formulate the result, we recall some definitions.

A linear operator $P\colon L^1(X)\to L^1(X)$ is said to be a \emph{Markov operator} if $Pf\geq 0$ and $\|Pf\|=\|f\|$ for every $f\in L^1(X)$ such that $f\geq 0$. It is easy to see that every Frobenius--Perron operator is a special type of Markov operator. We say that a sequence $(f_m)_{m\in\mathbb N}$ of functions from $L^1(X)$ is \emph{weakly Ces\`{a}ro convergent} to a function $f\in L^1(X)$ if
\begin{equation*}
\lim_{m\to\infty}\frac{1}{m}\sum_{k=1}^m\int_Xf_k(x)h(x)d\mu(x)=\int_Xf(x)h(x)d\mu(x)
\end{equation*}
for every $h\in L^{\infty}(X)$. A Markov operator $P\colon L^1(X)\to L^1(X)$ such that $P1=1$ is said to be \emph{ergodic} if for every density $f\in L^1(X)$ the sequence $(P^mf)_{m\in\mathbb N}$ is weakly Ces\`{a}ro convergent to $1$. %Since

Note that if $P\colon L^1(X)\to L^1(X)$ is a Markov operator and $\varphi\in L^1(X)$ is a solution of equation~\eqref{E}, then
\begin{equation*}
\int_Xg(x)d\mu(x)=\int_X\varphi(x)d\mu(x)-\int_XP\varphi(x)d\mu(x)=0.
\end{equation*}
Thus condition~\eqref{g} is necessary for $g\in L^1(X)$ in order that equation~\eqref{E} has a solution in $L^1(X)$.

\begin{theorem}\label{thm51}
Assume~\eqref{g} and let $P\colon L^1(X)\to L^1(X)$ be an ergodic Markov operator. Then equation~\eqref{E} has a solution in $L^1(X)$ if and only if the sequence $\big(\sum_{k=0}^{m-1}\frac{m-k}{m}P^{k}g\big)_{m\in\mathbb N}$ converges in $L^1(X)$. Moreover, every solution $\varphi\in L^1(X)$ of equation~\eqref{E} is of the form
\begin{equation*}
\varphi=\lim_{m\to\infty}\sum_{k=0}^{m-1}\frac{m-k}{m}P^{k}g+c,
\end{equation*}
where $c$ is a real constant.
\end{theorem}

\begin{pf}
Since $P$ is an ergodic Markov operator, it follows that $1$ is the unique density such that $P1=1$; indeed, assuming that there exists another density $f\in L^1(X)$ such that $Pf=f$, we would have
\begin{align*}
\int_Xf(x)h(x)d\mu(x)&=\lim_{m\to\infty}\frac{1}{m}\sum_{k=1}^m\int_XP^kf(x)h(x)d\mu(x)=\int_Xh(x)d\mu(x)
\end{align*} 	
for every $h\in L^\infty(X)$, which is impossible in the case where $f\neq 1$. Now from \cite[Theorem~5.2.2]{LM1994} (see also Proposition~5.2.1 in the same source) we conclude that for every density $f\in L^1(X)$ the sequence $\big(\frac{1}{m}\sum_{k=1}^{m}P^kf\big)_{m\in\mathbb N}$ converges to $1$ in $L^1(X)$, and by the linearity of $P$ we deduce that for every $f\in L^1(X)$ the sequence $\big(\frac{1}{m}\sum_{k=1}^{m}P^kf\big)_{m\in\mathbb N}$ converges to $\int_Xf(x)d\mu(x)$ in $L^1(X)$.
In particular, making use of~\eqref{g}, that is that the integral of $g$ over~$X$ vanishes, we obtain
\begin{equation}\label{Pg}
\lim_{m\to\infty}\frac{1}{m}\sum_{k=1}^{m}P^kg=0.
\end{equation}
	
Assume first that the sequence $\big(\sum_{k=0}^{m-1}\frac{m-k}{m}P^{k}g\big)_{m\in\mathbb N}$ converges in $L^1(X)$. Fix a real constant $c$ and set $\varphi=\lim_{m\to\infty}\sum_{k=0}^{m-1}\frac{m-k}{m}P^{k}g+c$. The linearity and continuity of $P$ jointly with  the equality $P1=1$ and \eqref{Pg} imply
\begin{align*}
P\varphi+g&=\lim_{m\to\infty}\sum_{k=0}^{m-1}\frac{m-k}{m}P^{k+1}g+Pc+g\\
&=\lim_{m\to\infty}\sum_{k=0}^{m-1}\frac{m-k}{m}P^{k}g+\lim_{m\to\infty}\frac{1}{m}\sum_{k=1}^{m}P^{k}g+c=\varphi.
\end{align*}
	
Assume now that $\varphi\in L^1(X)$ satisfies~\eqref{E}. Then, by the linearity of~$P$, we have
\begin{equation*}
\frac{m-k}{m}P^k\varphi=\frac{m-k}{m}P^{k+1}\varphi+\frac{m-k}{m}P^kg
\end{equation*}
for all $m\in\mathbb N$ and $k\in\{0,\ldots,m-1\}$. Adding the above equation over $k=0,\ldots,m-1$ with fixed $m$, leads to
\begin{eqnarray*}
\varphi-\frac{1}{m}\sum_{k=1}^{m}P^k\varphi=\sum_{k=0}^{m-1}\frac{m-k}{m}P^{k}g
\end{eqnarray*}
for every $m\in\mathbb N$. Finally, passing with $m$ to $\infty$ and making use of the fact that the sequence $\big(\frac{1}{m}\sum_{k=1}^{m}P^k\varphi\big)_{m\in\mathbb N}$ converges to $\int_X\varphi(x)d\mu(x)$ in $L^1(X)$,  we conclude that the sequence $\big(\sum_{k=0}^{m-1}\frac{m-k}{m}P^{k}g\big)_{m\in\mathbb N}$ converges in $L^1(X)$ and that
\begin{eqnarray*}
\varphi-\int_X\varphi(x)d\mu(x)=\lim_{m\to\infty}\sum_{k=0}^{m-1}\frac{m-k}{m}P^{k}g,
\end{eqnarray*}
which completes the proof. \hfill$\square$
\end{pf}

Theorem~\ref{thm51} generalizes Theorem~\ref{thm31} in two directions, because there are Markov operators that are not Frobenius--Perron operators and there are transformations $S$ that are not exact, but the corresponding Frobenius--Perron operators are ergodic. For example, it is easy to see that the operator defined in Section~3 by formula~\eqref{PS} can be ergodic, but not exact. Moreover, it fails to be a Frobenius--Perron operator in the case where at least one of the functions $f_n$ does not satisfy the Luzin's condition~N, but it is still a Markov operator in such a case.

The second result of this section shows that it can happen that equation~\eqref{E} can have exactly one integrable solution; note that in such a case we must have $P1\neq 1$.

%\begin{theorem}\label{thm52}
%Assume that $P\colon L^1(X)\to L^1(X)$ is a linear and continuous operator such that for every density $f\in L^1(X)$ the sequence $(P^mf)_{m\in\mathbb N}$ converges to the trivial function in $L^1(X)$. Then equation~\eqref{E} has a solution in $L^1(X)$ if and only if the series $\sum_{m=0}^{\infty}P^m g$ converges in $L^1(X)$. Moreover, every solution $\varphi\in L^1(X)$ of equation~\eqref{E} is of the form
%\begin{equation}\label{varphi}
%\varphi=\sum_{m=0}^{\infty}P^m g.
%\end{equation}
%\end{theorem}

\begin{theorem}\label{thm52}
Assume that the operator $P\colon L^1(X)\to L^1(X)$ is linear and continuous such that for every density $f\in L^1(X)$ the sequence $(P^mf)_{m\in\mathbb N}$ converges to the trivial function in $L^1(X)$. Then equation~\eqref{E} has a solution in the space~$L^1(X)$ if and only if the series $\sum_{m=0}^{\infty}P^m g$ converges in $L^1(X)$. Moreover, every solution $\varphi\in L^1(X)$ of equation~\eqref{E} is of the form
\begin{equation}\label{varphi}
\varphi=\sum_{m=0}^{\infty}P^m g.
\end{equation}
\end{theorem}

\begin{pf}	
By the linearity of $P$ it is easy to see that the sequence $(P^mf)_{m\in\mathbb N}$ converges to the trivial function in $L^1(X)$ for every $f\in L^1(X)$.
	
Assume first that the series $\sum_{m=0}^{\infty}P_S^m g$ converges in $L^1(X)$. Setting $\varphi=\sum_{m=0}^{\infty}P^m g$ and applying the linearity and continuity of $P$ we obtain
\begin{equation*}
P\varphi+g=\sum_{m=0}^{\infty}P^{m+1}g+g=\sum_{m=0}^{\infty}P^m g=\varphi.
\end{equation*}
	
Assume now that $\varphi\in L^1(X)$ satisfies~\eqref{E}. By the linearity of $P$ we have $P^k\varphi=P^{k+1}\varphi+P^kg$ and hence
\begin{equation*}
\sum_{k=0}^m P^k g=\varphi-P^{m+1}\varphi
\end{equation*}
for every $m\in\mathbb N$. Passing with $m$ to $\infty$ we deduce that the series $\sum_{k=0}^\infty P^k g$ converges in $L^1(X)$ and that \eqref{varphi} holds.\hfill$\square$
\end{pf}

To give an example of a realization of the assumptions of Theorem~\ref{thm52} fix, to the end of this paper, strictly monotone functions $f_1,\ldots,f_N\colon[0,1]\to[0,1]$ satisfying condition~\eqref{c2} and consider the operator $P_0\colon L^1([0,1])\to L^1([0,1])$ defined by
\begin{equation*}
P_0f=\sum_{n=1}^{N}|f_n'|\big(f\circ f_n\big).
\end{equation*}
Obviously, $P_0$ is linear. To see that $P_0$ is continuous note that \eqref{c2} yields
\begin{equation}\label{p}
\int_AP_0f(x)dx=\sum_{n=1}^{N}\int_{A}|f_n'(x)|f(f_n(x))dx=\int_{\bigcup_{n=1}^{N}f_n(A)}f(y)dy
\end{equation}
for all nonnegative $f\in L^1([0,1])$ and Lebesgue measurable sets $A\subset[0,1]$.

Assume now that the family $\{f_1,\ldots,f_N\}$ forms an iterated function system and let $A_*$ be its \emph{attractor}, i.e.
\begin{equation*}
A_*=\bigcap_{m\in\mathbb N}A_m,
\end{equation*}
where $A_0=[0,1]$ and $A_{m}=\bigcup_{n=1}^{N}f_n(A_{m-1})$ for every $m\in\mathbb N$. Fix a nonnegative $f\in L^1([0,1])$. According to~\eqref{p} we have
\begin{equation*}
\|P_0^mf\|=\int_{A_m}f(y)dy
\end{equation*}
for every $m\in\mathbb N$, and as the sequence $(A_m)_{m\in\mathbb N}$ is descending we get
\begin{equation*}
\lim_{m\to\infty}\|P_0^mf\|=\int_{A_*}f(y)dy.
\end{equation*}
In consequence, we have proved the following lemma.

\begin{lemma}\label{lem53}
If the family $\{f_1,\ldots,f_N\}$ forms an iterated function system with the attractor of Lebesgue measure zero, then for every nonnegative $f\in L^1([0,1])$ the sequence  $(P_0^mf)_{m\in\mathbb N}$ converges to the trivial function in $L^1([0,1])$.
\end{lemma}

As an immediate consequence of Theorem \ref{thm52} and Lemma~\ref{lem53}, we obtain the following result.

\begin{corollary}\label{cor54}
Assume that the family $\{f_1,\ldots,f_N\}$ forms an iterated function system and let its attractor have Lebesgue measure zero. Then equation~\eqref{e} has a solution in $L^1([0,1])$ if and only if the series $\sum_{m=0}^{\infty}P_0^m g$ converges in $L^1([0,1])$. Moreover, every solution $\varphi\in L^1([0,1])$ of equation~\eqref{e} is of the form
\begin{equation*}
\varphi=\sum_{m=0}^{\infty}P_0^m g.
\end{equation*}
\end{corollary}

Observe that assumptions of Corollary~\ref{cor54} are satisfied if the $f_n$ are defined by~\eqref{affine} with real numbers $\alpha_1,\ldots,\alpha_N$ and non-zero real numbers $\beta_1,\ldots,\beta_N$ such that
\begin{align*}\label{1a}
0&\leq\min\{\beta_1,\alpha_1+\beta_1\}<\max\{\beta_1,\alpha_1+\beta_1\}\leq\min\{\beta_2,\alpha_2+\beta_2\}\\
&<\max\{\beta_2,\alpha_2+\beta_2\}\leq\cdots\leq\min\{\beta_N,\alpha_N+\beta_N\}\\
&<\max\{\beta_N,\alpha_N+\beta_N\}\leq 1
\end{align*}
and
\begin{equation*}\label{1b}
\bigcup_{n=1}^{N}[\min\{\beta_n,\alpha_{n}+\beta_{n}\},\max\{\beta_n,\alpha_{n}+\beta_{n}\}]\neq[0,1].
\end{equation*}
It is easy to see that the family  $\{f_0,\ldots,f_N\}$ forms an iterated function system and its attractor $A_*$ has Lebesgue measure zero. Thus by Corollary~\ref{cor54} we conclude (in contrast to the counterpart case from Section 3) that now equation~\eqref{e1} has exactly one solution $\varphi\in L^1([0,1])$ and it is of the form~\eqref{sol} with $c=0$.

\section*{Acknowledgements}
This research was supported by the University of Silesia Mathematics Department (Iterative Functional Equations and Real Analysis program).
%\section*{References}

\bibliographystyle{plain}
\bibliography{Invariant}

\begin{thebibliography}{10}

\bibitem{Dud2002}
R.~M. Dudley.
\newblock {\em Real analysis and probability}, volume~74 of {\em Cambridge
  Studies in Advanced Mathematics}.
\newblock Cambridge University Press, Cambridge, 2002.
\newblock Revised reprint of the 1989 original.

\bibitem{LM1994}
Andrzej Lasota and Michael~C. Mackey.
\newblock {\em Chaos, fractals, and noise}, volume~97 of {\em Applied
  Mathematical Sciences}.
\newblock Springer-Verlag, New York, second edition, 1994.
\newblock Stochastic aspects of dynamics.

\bibitem{M1985}
Janusz Matkowski.
\newblock Remark on {BV}-solutions of a functional equation connected with
  invariant measures.
\newblock {\em Aequationes Math.}, 29(2-3):210--213, 1985.

\bibitem{MW2012}
Jolanta Misiewicz and Jacek Weso{\l}owski.
\newblock Winding planar probabilities.
\newblock {\em Metrika}, 75(4):507--519, 2012.

\bibitem{MZ2018}
Janusz Morawiec and Thomas Z\"{u}rcher.
\newblock On a problem of {J}anusz {M}atkowski and {J}acek {W}eso\l owski.
\newblock {\em Aequationes Math.}, 92(4):601--615, 2018.

\bibitem{MZ}
Janusz Morawiec and Thomas Z\"urcher.
\newblock On a problem of {Janusz} {Matkowski} and {Jacek} {Weso\l
  owski},~{II}, 2018.
\newblock \url{https://arxiv.org/abs/1801.07944}.

\bibitem{N1991}
Kazimierz Nikodem.
\newblock On {$\epsilon$}-invariant measures and a functional equation.
\newblock {\em Czechoslovak Math. J.}, 41(116)(4):565--569, 1991.

\bibitem{P1977}
K.~R. Parthasarathy.
\newblock {\em Introduction to probability and measure}.
\newblock The Macmillan Co. of India, Ltd., Delhi, 1977.

\bibitem{R1957}
A.~R{\'e}nyi.
\newblock Representations for real numbers and their ergodic properties.
\newblock {\em Acta Math. Acad. Sci. Hungar}, 8:477--493, 1957.

\bibitem{Ro1961}
V.~A. Rohlin.
\newblock Exact endomorphisms of a {L}ebesgue space.
\newblock {\em Izv. Akad. Nauk SSSR Ser. Mat.}, 25:499--530, 1961.

\bibitem{rohlin1964exact}
Vladimir~A Rohlin.
\newblock Exact endomorphisms of a {L}ebesgue space.
\newblock {\em Amer. Math. Soc. Transl. Ser}, 2(39):1--36, 1964.

\end{thebibliography}

\end{document}